\documentclass[final, 1p]{elsarticle}

\usepackage{amsmath}
\usepackage{amsfonts}
\usepackage{amssymb}
\usepackage{mathrsfs}	
\usepackage{bm}
\usepackage{fancybox}
\usepackage{makeidx}
\usepackage{graphicx}
\usepackage{subfigure}
\usepackage{epsfig}
\usepackage{verbatim}
\usepackage{hyperref}
\usepackage{cases}
\usepackage{booktabs}
\usepackage{upgreek}
\usepackage{tabularx,multirow}
\usepackage{color}
\usepackage[active]{srcltx}
\usepackage{citesort}
\usepackage{ulem}   
\let \oldmarginpar \marginpar
\renewcommand{\marginpar}[1]{\oldmarginpar{\color{red}{#1}}}

\newcounter{algleo}
\newlength{\lefttab}
\newlength{\numberoffset}
\newenvironment{algleo}%
  {\trivlist
   \topsep=0pt\parsep=0pt\itemsep=0pt
   \def\li{\item\refstepcounter{algleo}\makebox[0.8em][r]{\thealgleo\hspace{\numberoffset}}
       \hangafter1\hangindent1.8em\noindent}%
   \def\linonumber{\item\makebox[0.8em][r]{\hspace{\numberoffset}}
       \hangafter1\hangindent1.8em\noindent}%
   \addtolength{\lefttab}{1.25em}
   \addtolength{\numberoffset}{1.25em}
   \leftskip=\lefttab}%
  {\endtrivlist}

\definecolor{MyGray}{rgb}{0.5, 0.5, 0.5}
\definecolor{MyDarkGreen}{rgb}{0.0, 0.46, 0.29}
\definecolor{MyLightBlue}{rgb}{0.85, 0.85, 1.0}
\definecolor{MyLightGreen}{rgb}{0.9, 1.0, 0.9}

\begin{document}

\begin{frontmatter}

\title{A fast directional boundary element method for
	high frequency acoustic problems in three dimensions}

\author[]{Yanchuang Cao}
\ead{caoyanch@126.com}

\author[]{Lihua Wen\corref{cor1}}
\ead{lhwen@nwpu.edu.cn}

\author[]{Jinyou Xiao}
\ead{xiaojy@nwpu.edu.cn}

\cortext[cor1]{Corresponding author}
\address{College of Astronautics, Northwestern Polytechnical University, 
	Xi'an 710072, P. R. China}

\begin{abstract}
A highly efficient fast boundary element method (BEM) for
solving large-scale engineering acoustic problems in a broad frequency 
range is developed and implemented. The acoustic problems are modeled by 
the Burton-Miller boundary integral equation (BIE), thus the fictitious 
frequency issue is completely avoided. The BIE is discretized by using 
the collocation method with piecewise constant elements. The linear 
systems are solved iteratively 
and accelerated by using a newly developed kernel-independent wideband 
fast directional algorithm (FDA) for fast summation of oscillatory 
kernels. In addition, the computational efficiency of the FDA is further 
promoted by exploiting the low-rank features of the translation matrices. 
The high accuracy and nearly linear 
computational complexity of the present method are clearly
demonstrated by typical examples. An acoustic scattering problem with
dimensionless wave number $kD$ (where $k$ is the wave number and $D$ is 
the typical length of the obstacle) up to 1000 and the degrees of freedom 
up to 4 million is successfully solved within 4 hours on a computer with 
one core and the memory usage is 24.7 GB.
\end{abstract}


\end{frontmatter}

\section{Introduction}

Acoustic wave propagation is a commonly studied problem for its wide application in
noise controls, ultrasonic diagnostics, sonar imaging, etc. One of the most popular
solving approaches is boundary element method due to its unique advantages,
such as dimension reduction, high accuracy and suitability for infinite domain
cases. The conventional boundary element method leads to dense system matrix.
As a result, the computational complexity is at least of order $O(N^2)$, which makes
it prohibitive for large-scale problems.

In the past three decades, many fast algorithms have been proposed to circumvent this
disadvantage. A group of these algorithms make use of the asymptotic smooth 
property of the kernel function when the field points are far away from the source 
points. As a result, the large submatrices of the system matrix corresponding to 
far-field interactions are of numerically low rank. These algorithms are 
performed on an octree (quadtree for two dimensional cases), which is 
constructed by subdividing the integral region recursively. By applying 
low rank approximations hierarchically on the octree, the computational 
complexity can be reduced to $O(N)$. The low rank decomposition of the 
submatrices can be generated by various methods, resulting in 
diverse algorithms, including fast multipole method (FMM)\cite{fmm, pwfmm, %
fmm_Darve2000}, $\mathcal{H}$-matrix method\cite{Hmatrix}, panel clustering 
method\cite{pnlclst}, ACA\cite{aca, aca_BR2003}, wavelet compression 
method\cite{wbem, twbem} etc. However, it is found that these low-rank 
approximating algorithms are only suitable for low frequency problems, since 
the ranks of the submatrices tend to be proportional to their sizes, %
leading to $O(N^2)$ complexity for high frequency problems.

Another class of these fast algorithms make use of the translational
invariant property of the kernel function. By mapping the information on the 
elements onto Cartesian grids, then diagonalizing the matrix by Fast Fourier 
Transformation\cite{pfft, pfft_YanzyGaoxw}, the required operations 
can be reduced significantly. It can be very efficient for
both low and high frequency problems, but its computational complexity is 
$O(N^{4/3} \log N)$ when applied in accelerating 
BEM\cite{pfft_Bruno2001, pfft_Bruno2012}. Besides, it is well known
that the performance of the FFT-based algorithms deteriorates 
in the cases with highly nonuniform mesh discretization.

Using the diagonal forms of the translation operators in FMM 
makes it possible to obtain $O(N \log N)$ complexity for high 
frequency problems\cite{diagonalform}. Later the wideband 
FMM\cite{wbfmm, wbfmm_Gumerov2009} is developed which 
successfully avoids its numerical instability at low frequencies 
by combining it with the traditional FMM. 
In wideband FMM, the octree is divided into high frequency regime and 
low frequency regime according to the size of the cubes in each level. 
Different translation methods are applied in these regimes, i.e., in the low 
frequency regime, the translations are performed in the same way as in 
traditional FMM; while far field signature and diagonal forms of the translation
operators are used in the high frequency regime. It is stable, accurate and 
efficient for both low and high frequency problems. However, it requires 
analytical expansions of the kernel, thus the translations are very 
complicated and the algorithm is kernel-dependent. This poses a severe 
limitation on its applications.

Fast directional algorithm is another efficient algorithm for solving high 
frequency problems\cite{fda, fda2010}, by which the computational complexity 
can also be reduced to $O(N \log N)$. It is a FMM-like algorithm which takes 
the advantage of the directional low rank property of the kernel function
to do the translations in the high frequency regime. Consequently, the low
rank approximation are applied in both the low and high frequency regimes, %
the only difference is the definition of the interaction list. Various 
low rank approximating techniques can be used, resulting in variants 
of the fast directional algorithm\cite{fda2010, fda_Chebyshev, dirH2}. 

In this paper, the fast directional algorithm based on equivalent 
densities is adapted to accelerate the 3D acoustic BEM for the first
time, and a modified version of a recently developed M2L translations
accelerating technique named as SArcmp is applied to improve its efficiency. 
The advantages of this algorithm lies in the following two aspects. 
First, no analytical expansions for the kernel function is 
required, thus the algorithm is completely kernel-independent and 
easy to implement. Second, integrals with different layer kernel 
functions can be accelerated by the same process, thus it is 
very convenient to handle Burton-Miller formulation
in which four layer kernel functions are included. In this sense, it 
is more suitable to accelerated acoustic problems than other 
fast algorithm.

\section{Boundary integral formulation for acoustic problems}

Consider the acoustic problems described by Helmholtz equation
\begin{equation}
	\nabla^2 u(\bm{x}) + k^2 u(\bm{x}) = 0, \quad \bm{x} \in \Omega,
\end{equation}
where $u$ is the velocity potential, $k = \omega/c$ is the wavenumber, and 
$\Omega$ is the acoustic field domain. 
The acoustic field can be solved by conventional boundary integral 
equation
\begin{equation} \label{cbie}
	c(\bm{x}) u(\bm{x}) + \int_\Gamma \frac{\partial G(\bm{x}, \bm{y})}
	{\partial \bm{n}(\bm{y})} u(\bm{y}) \mathrm{d} \bm{y} = \int_\Gamma 
	G(\bm{x}, \bm{y}) \frac{\partial u(\bm{y})}{\partial \bm{n}(\bm{y})} 
	\mathrm{d} \bm{y} + u^\text{inc}(\bm{x}), \quad \bm{x} \in \Gamma,
\end{equation}
where $\Gamma = \partial \Omega$ is the boundary of the acoustic 
field, $c(\bm{x})$ is the solid angle at $\bm{x}$, and $G(\bm{x}, \bm{y})$ is
the fundamental solution
\begin{equation}
	G(\bm{x}, \bm{y}) = \frac{e^{ikr}}{4\pi r}, \quad r = | \bm{x} - \bm{y} |.
\end{equation}
However, it is well known that it fails to yield unique solutions
at the characteristic frequencies. One of the most widely used methods 
to overcome this problem is the Burton-Miller formulation, which
is generated by combining \eqref{cbie} and its normal derivatives 
\begin{equation} \label{eq_bm}
\begin{split}
	c(\bm{x}) u(\bm{x}) + \int_\Gamma \frac{\partial G(\bm{x}, \bm{y})}
	{\partial \bm{n}(\bm{y})} u(\bm{y}) \mathrm{d} \bm{y} + 
	\alpha \int_\Gamma \frac{\partial^2 G(\bm{x}, \bm{y})} 
	{\partial \bm{n}(\bm{x}) \partial \bm{n}(\bm{y})} 
	u(\bm{y}) \mathrm{d} \bm{y} \\
	= -\alpha c(\bm{x}) \frac{\partial u(\bm{x})}{\partial \bm{n}(\bm{x})} 
	+ \int_\Gamma G(\bm{x}, \bm{y}) \frac{\partial u(\bm{y})}
	{\partial \bm{n}(\bm{y})} \mathrm{d} \bm{y} + 
	\alpha \int_\Gamma \frac{\partial G(\bm{x}, \bm{y})}{\partial 
		\bm{n}(\bm{x})} \frac{\partial u(\bm{y})}
	{\partial \bm{n}(\bm{y})} \mathrm{d} \bm{y} \\
	+ u^\text{inc}(\bm{x}) + \alpha \frac{\partial u^\text{inc}(\bm{x})}
	{\partial \bm{n}(\bm{x})}, \quad \bm{x} \in \Gamma,
\end{split}
\end{equation}
where $\alpha$ is the combining factor that is suggested to be chosen as 
$i/k$\cite{condBM}. 

By discretizing \eqref{eq_bm} with basis functions $\chi(\bm{x})$ and weight
functions $w(\bm{x})$, the integrals would be transformed into summations.
Take the left hand side for example, it can be discretized into 
\begin{equation} \label{eq_sum}
\begin{split}
	p_i = & \sum_{j=1}^N \int_{\Gamma_i} w_i(\bm{x}) \int_{\Gamma_j} \left[ 
	\frac{\partial G(\bm{x}, \bm{y})}{\partial \bm{n}(\bm{y})} 
	+ \alpha \frac{\partial^2 G(\bm{x}, \bm{y})}{\partial \bm{n}(\bm{x}_i) 
	\partial \bm{n}(\bm{y})} \right]	q_j \chi_j(\bm{y}) \mathrm{d} \bm{y} 
	\mathrm{d} \bm{x}\\
	= & \int_{\Gamma_i} w_i(\bm{x}) \left[ 1 + \alpha \frac{\partial}{\partial 
		\bm{n}(\bm{x})} \right] \left[ \sum_{j=1}^N \int_{\Gamma_j}
	\frac{\partial G(\bm{x}, \bm{y})}{\partial \bm{n}(\bm{y})} 
	\chi_j(\bm{y}) \mathrm{d} \bm{y} \cdot q_j \right] \mathrm{d} \bm{x},\\ 
	&\quad i = 1, 2, \dots, N,
\end{split}
\end{equation}
where $\Gamma_i$ and $\Gamma_j$ is the supporting region of the $i$-th 
weight function $w_i(\bm{x})$ and the $j$-th basis function 
$\chi_j(\bm{y})$, respectively, and $q_j$ is the coefficient for 
$\chi_j(\bm{y})$. Evaluating the summation directly requires $O(N^2)$ 
operators. In the next section, we will discuss how to 
accelerate the evaluation by fast directional algorithm.

\section{Fast directional algorithm for Burton-Miller formulation}

The key of fast directional algorithm is the construction of 
the fast potential evaluating scheme using the directional low rank 
property of the kernel function. In this section, the fast evaluating
scheme is also based on equivalent densities and check potentials as
\cite{fda, fda2010}. However, the equivalent points and check points are 
defined as the quadrature points instead of defined by pseudo skeleton 
approach, resulting in a simpler fast directional algorithm.

\subsection{Directional low rank approximation}

\begin{figure}[h]
	\centering
	\includegraphics[width=0.6\textwidth]{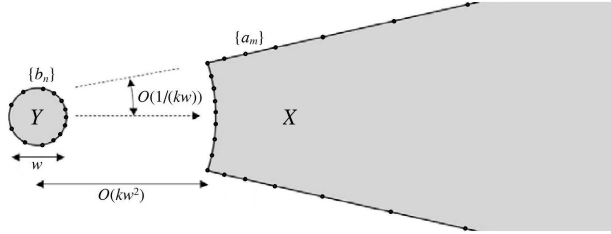}
	\caption{Source set $Y$ and target set $X$ satisfying the directional
	parabolic separation condition.}
	\label{fig_para_sep}
\end{figure}

Suppose $X$ and $Y$ be the target point set and the source point set, respectively.
When $X$ and $Y$ satisfy the \emph{directional parabolic separation 
condition}, as shown in Figure \ref{fig_para_sep}, the kernel function can
be approximated by
\begin{equation} \label{eq_dirknlrep}
	\left| G(\bm{x}, \bm{y}) - \sum_{i=1}^{T(\varepsilon)} \alpha_i(\bm{x}) 
	\beta_i(\bm{y}) \right| < \varepsilon,
\end{equation}
where $T(\varepsilon)$ has an upper bound that is independent of $k$ and $w$.
In this case, the evaluation for the potentials 
on $X$ can be accelerated via equivalent densities and check 
potentials. 

First we find certain distributions of monopole sources $\sigma(\bm{y})$ on 
\emph{directional outgoing equivalent surface} $\bm{y}^{Y, \text{o}}$, that 
can reproduce the potential field $p(\bm{x})$ in $X$ excited by
arbitrary source densities $q(\bm{y})$ in $Y$.
Similar with the kernel independent FMM\cite{kifmm}, the equivalent 
surface $\bm{y}^{Y, \text{o}}$ need to enclose $Y$ in order to guarantee 
the existence of the \emph{directional outgoing equivalent densities} 
$\sigma(\bm{y})$; and a check surface $\bm{x}^{Y, \text{o}}$ enclosing 
$X$ can be defined, such that once the $\sigma(\bm{y})$
can reproduce the potential field on $\bm{x}^{Y, \text{o}}$, it can 
excite the same potential field inside $X$. Therefore, the 
\emph{directional outgoing equivalent densities}
$\sigma(\bm{y})$ on $\bm{y}^{Y, \text{o}}$ can be computed by 
the following equation
\begin{equation} \label{eq_eqvdencmp}
	\int_{\bm{y}^{Y, \text{o}}} G(\bm{x}, \bm{y}) \sigma(\bm{y}) \mathrm{d} \bm{y}
	= p^{Y, \text{o}}(\bm{x}), \quad \bm{x} \in \bm{x}^{y, \text{o}}.
\end{equation}
where $p^{Y, \text{o}}(\bm{x})$ is calculated by the original source 
densities $q(\bm{y})$ in $Y$. For point sources, $p^{Y, \text{o}}(\bm{x})$ 
is evaluated by summation; while for 
distributed sources, it should be evaluated by quadrature. Equation 
\eqref{eq_eqvdencmp} can be viewed as a transformation from 
$\sigma(\bm{y})$ to $p^{Y, \text{o}}(\bm{x})$, and 
the transformation is of rank $T(\varepsilon)$. Therefore, it can be discretized
by Nystr\"{o}m method using $O(T(\varepsilon))$ \emph{directional outgoing 
equivalent points} $\bm{y}$'s and $O(T(\varepsilon))$ \emph{directional outgoing 
check points} $\bm{x}$'s, i.e., 
\begin{equation}
	p^{Y, \text{o}}(\bm{x}_i) = \int_Y G(\bm{x}_i, \bm{y}) \sigma(\bm{y}) 
	\mathrm{d} \bm{y}
	= \sum_j w_j G(\bm{x}_i, \bm{y}_j) \sigma(\bm{y}_j), \quad
	\bm{x}_i \in \bm{x}^{Y, \text{o}}, \bm{y}_j \in \bm{y}^{Y, \text{o}}.
\end{equation}
This suggests that we can take distribute monopole source densities 
$q^{Y, \text{o}}(\bm{y}_j) = w_j \sigma(\bm{y}_j)$ at 
quadrature points as the \emph{directional outgoing equivalent densities}, and
the transformation from \emph{directional outgoing equivalent densities} to
\emph{directional outgoing check potentials} becomes
\begin{equation} \label{eq_diseqvden}
	p^{Y, \text{o}}(\bm{x}_i) = \sum_j G(\bm{x}_i, \bm{y}_j) 
	q^{Y, \text{o}}(\bm{y}_j), \quad 
	\bm{x}_i \in \bm{x}^{Y, \text{o}}, \bm{y}_j \in \bm{y}^{Y, \text{o}}.
\end{equation}

\begin{figure}[h]
	\centering
	\subfigure[S2M]{
		\includegraphics[width=0.3\textwidth]{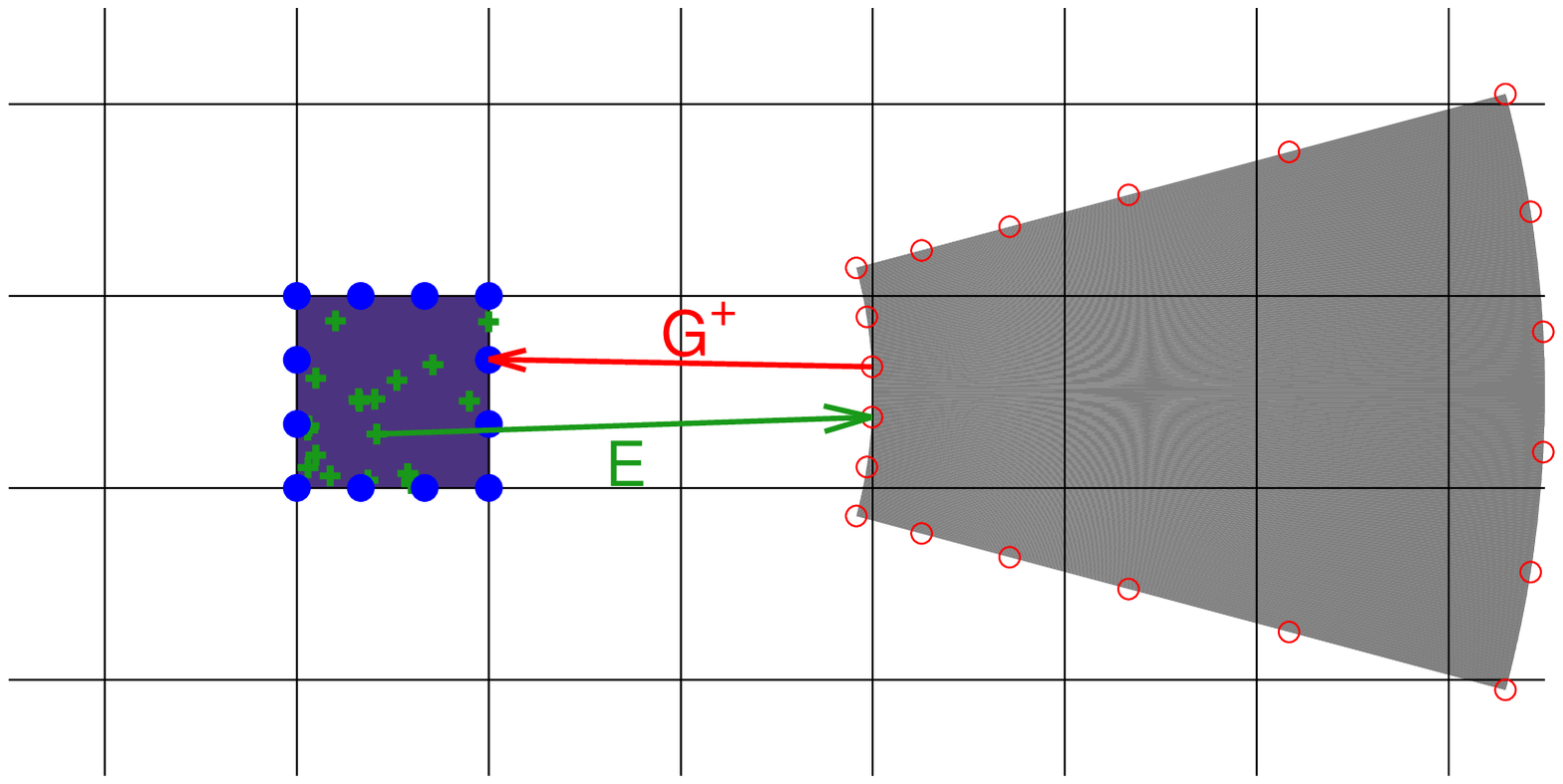}
		\label{fig_s2m}
	}
	\subfigure[M2L]{
		\includegraphics[width=0.3\textwidth]{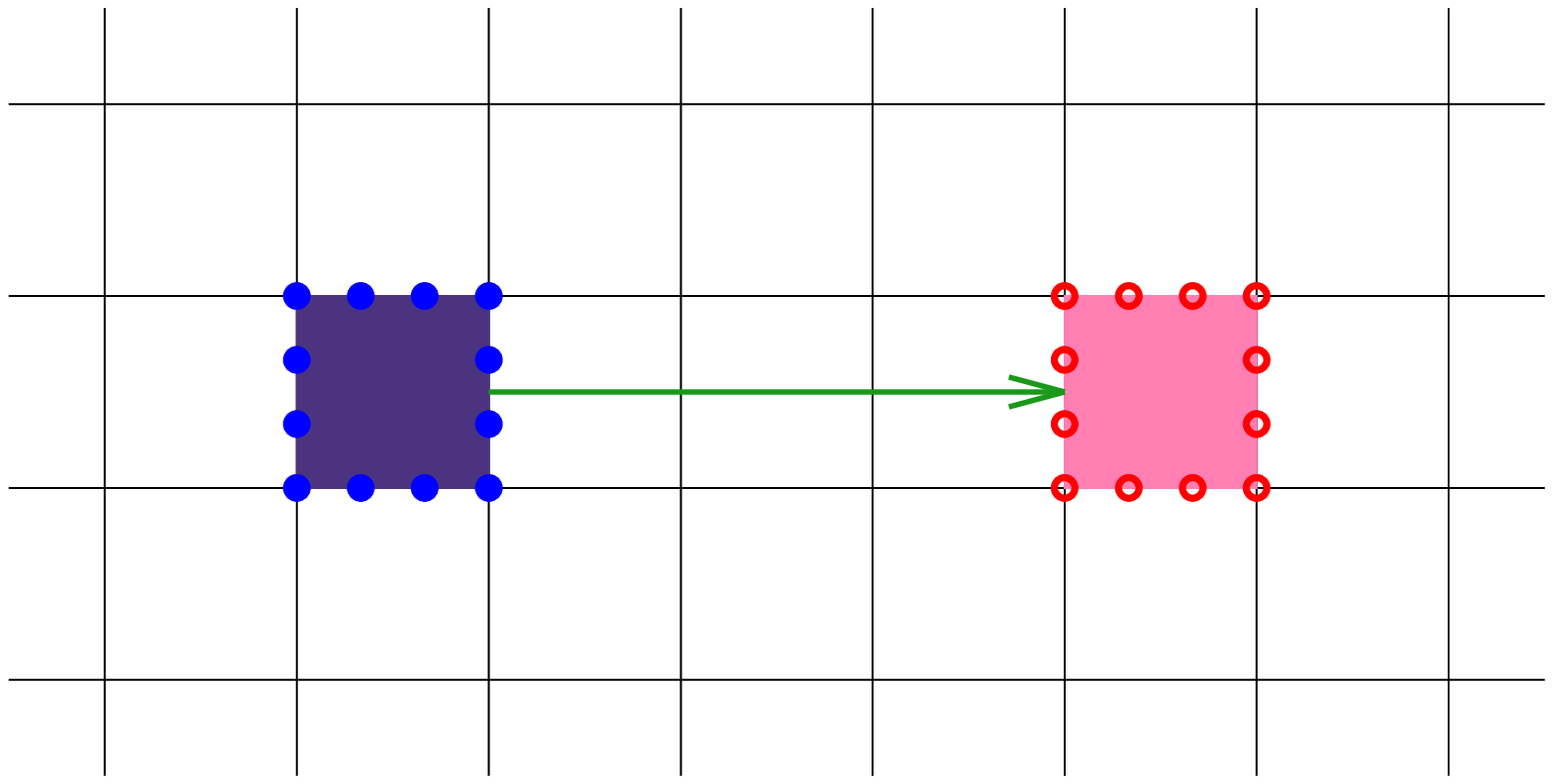}
		\label{fig_m2l}
	}
	\subfigure[L2T]{
		\includegraphics[width=0.3\textwidth]{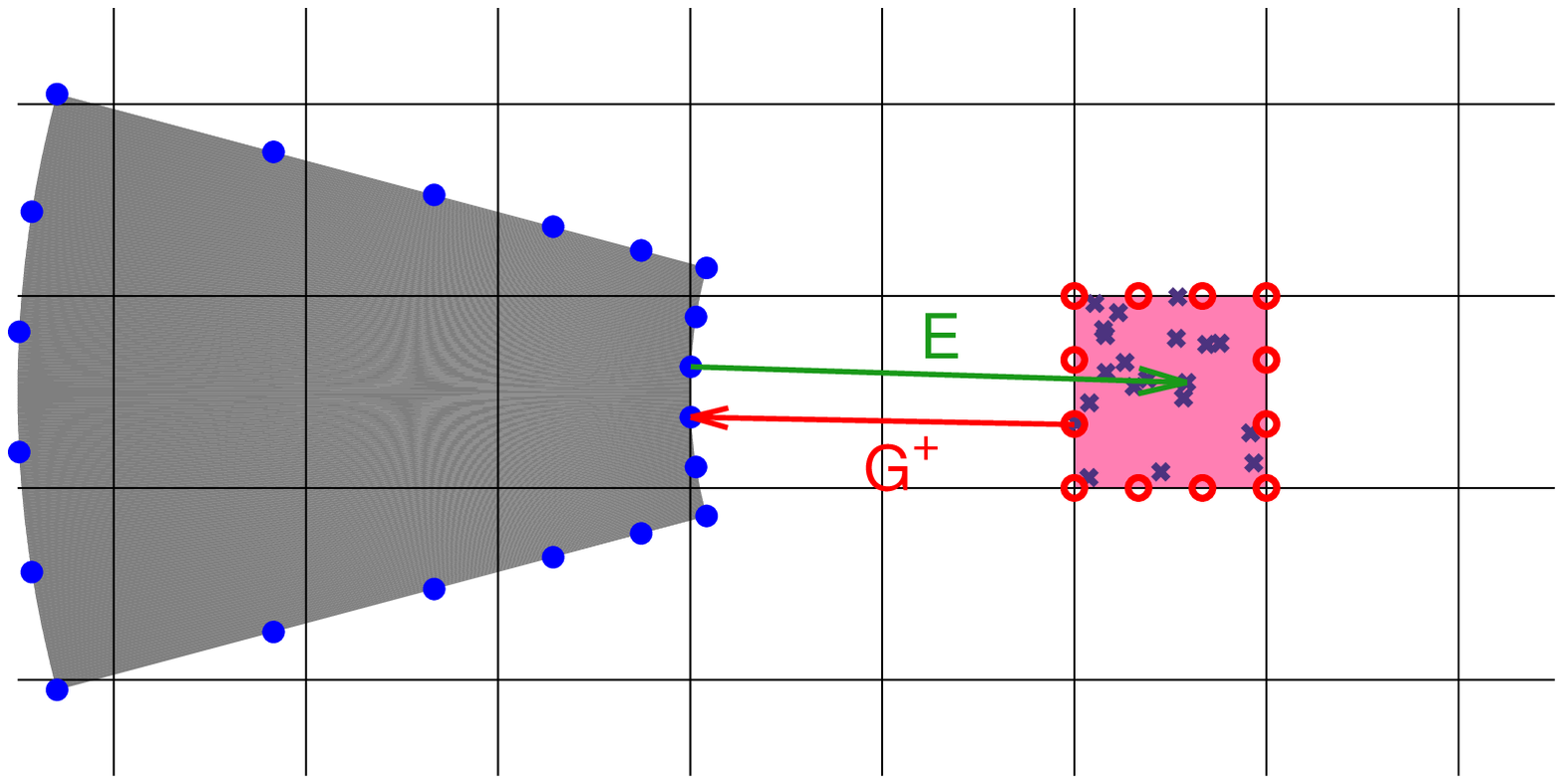}
		\label{fig_l2t}
	}
	\caption{Single level fast directional summation for $N$-body problems.}
	\label{fig_slfda}
\end{figure}

In our algorithm, the \emph{directional outgoing equivalent points} are distributed 
in the same way as the \emph{non-directional outgoing equivalent points} in 
\cite{kifmm, kifmbem}, that is, they are distributed on a cube surface with
$p$ points in each direction. The \emph{directional outgoing check points} are 
defined by mapping the points onto the surface of the directional cone
which is bounded by the size of the boundary $\Gamma$, and are focused at 
the smaller end, as illustrated in Figure \ref{fig_s2m}.

When the roles of $X$ and $Y$ are reversed, i.e., the potentials in $Y$ 
produced by sources inside $X$ need to be evaluated, the \emph{directional incoming 
equivalent densities} $q^{X, \text{i}}(\bm{y})$ on the \emph{directional 
incoming equivalent points} $\bm{y}^{X, \text{i}}$ can be constructed in 
the same manner via the \emph{directional incoming check potentials} 
$p^{X, \text{i}}(\bm{x})$ on the \emph{directional 
incoming check points} $\bm{x}^{X, \text{i}}$.

Since the kernel is rotational invariant, the \emph{directional
equivalent points} and the \emph{directional check points} for other
directional cones can be obtained by rotation, and the matrices 
translating the check potentials into equivalent densities remains 
the same.

Following the above scheme, the equivalent points and check 
points can be distributed straightforwardly instead of by 
pseudo skeleton approach as in \cite{fda2010}.
The transformations in high frequency regime can be accelerated in the same way as
that in low frequency regime, which has been discussed in detail in 
\cite{kifmm, kifmbem}, except that the outgoing check points and incoming 
equivalent points are \emph{directional}.

\subsection{Fast directional algorithm}

Similar with other fast directional algorithms, the octree is constructed 
by separating the computing domain recursively, and is divided into high 
frequency regime and low frequency regime. 
In our algorithm, the octree is constructed exactly in the same way as 
in FMM and kernel-independent FMM, thus there may be leaf cubes in the
high frequency regime. Therefore, our algorithm is completely \emph{adaptive},
while other fast directional algorithms\cite{fda, fda2010, fda_Chebyshev} 
are not since in those algorithms no adaptivity is used in the 
high frequency regime.

In the low frequency regime, the translations are computed by kernel-independent 
FMM\cite{kifmm, kifmbem}. In the high frequency, the interaction field of 
each cube is divided into directional cones, and the translations are 
accelerated by directional low rank approximation. 
Consider two high frequency leaf cubes $C$ and $D$ which are in the
same level and in each other's interaction list. 
We need to evaluate the potentials on $D$ generated by the sources 
in $C$. The accelerating approach is similar with the single level fast
directional algorithm for $N$-body problems, which is shown in Figure 
\ref{fig_slfda}, i.e., the evaluation can be accelerated by splitting 
it into three steps:
\begin{enumerate}
	\item Directional S2M translation: Compute the \emph{directional 
	outgoing equivalent densities}. First evaluate the \emph{directional outgoing
	check potentials} produced by the original source densities in cube $C$:
	\begin{equation} \label{eq_s2m}
		p^{C, \text{o}}(\bm{x}_i) = \sum_{j\in \sharp \Gamma_C} \int_{\Gamma_j} 
		\frac{\partial G(\bm{x}_i, \bm{y})}{\partial \bm{n}(\bm{y})} \chi_j(\bm{y})
		\mathrm{d}\bm{y} \cdot q_j, \quad \bm{x}_i \in \bm{x}^{C, \text{o}},
	\end{equation}
	where $\sharp \Gamma_C$ is the indices of the basis functions ``belonging
	to'' $C$. Then the \emph{directional outgoing equivalent densities}
	$q^{Y, \text{o}}(\bm{y}_j)$ can be computed by inverting \eqref{eq_diseqvden}.
	\item Directional M2L translation: Compute the \emph{directional 
	incoming check potentials}
	\begin{equation} \label{eq_m2l}
		p^{D, \text{i}}(\bm{x}_i) = \sum_j G(\bm{x}_i, \bm{y}_j) 
		q^{C, \text{o}}(\bm{y}_j), \quad
		\bm{x}_i \in \bm{x}^{D, \text{i}}, \bm{y}_j \in \bm{y}^{C, \text{o}}.
	\end{equation}
	\item Directional L2T translation: Compute the potentials on the $D$. First
	the \emph{directional incoming equivalent densities} 
	$q^{D, \text{i}}(\bm{y})$ are constructed by inversion, then evaluate the 
	potentials by the following equation
	\begin{equation} \label{eq_l2t}
	\begin{split}
		p_i &= \int_{\Gamma_i} w_i(\bm{x}) \left[ 1 + \alpha 
		\frac{\partial}{\partial \bm{n}(\bm{x})} \right] 
		\left[ \sum_j G(\bm{x}, \bm{y}_j) q^{D, \text{i}}(\bm{y}_j) 
		\right] \mathrm{d} \bm{x} \\
		&= \sum_j \int_{\Gamma_i} w_i(\bm{x}) \left[ G(\bm{x}, \bm{y}_j) +
		\alpha \frac{\partial G(\bm{x}, \bm{y}_j)}{\partial \bm{n}(\bm{x})} 
		\right] \mathrm{d} \bm{x} \cdot q^{D, \text{i}}(\bm{y}_j), 
		\quad \bm{y}_j \in \bm{y}^{D, \text{i}}.
	\end{split}
	\end{equation}
\end{enumerate}

In the multilevel fast directional algorithm, the M2M translation is
similar with S2M, but the \emph{directional outgoing check potentials} in $C$ 
are evaluated by the \emph{directional outgoing equivalent densities}
of $C$'s child cubes. The L2L translation is similar with L2T, but
instead of the potentials on the target points, the \emph{directional 
incoming check potentials} of $D$'s child cubes are evaluated. Thus, they 
are similar with that in KIFMM \cite{kifmm}, but should be transformed
into \emph{directional} in high frequency regime by using \emph{directional 
outgoing check points} and \emph{directional incoming equivalent points}.

Note that although there are two integrals in the left hand side of the 
Burton-Miller formulation, by using the above fast directional 
algorithm, Eq. \eqref{eq_sum} can be evaluated by one fast directional approach. 
Similarly, the evaluation of the integral in the right hand side can
also be accelerated by one fast directional approach.
Consequently, although there are four integrals in the Burton-Miller 
formulation, only two fast directional approaches are required. 

\subsection{Algorithm summary}

The overall multilevel fast directional algorithm accelerating \eqref{eq_sum}
is summarized as follows.

\setcounter{algleo}{0}
\begin{algleo}
    \linonumber \textbf{Algorithm} \quad 
	\textsc{Fast directional algorithm for \ref{eq_sum}}

    \linonumber \textsc{Step 1 Setup}
    \begin{algleo}
        \li Construct the octree adaptively.
        \li Define the far fields on each level. 
		\li Divide the octree into low and high frequency regimes. 
		\li Divide the far fields in the high frequency regime into 
		directional cones. 
		\li Construct interacting lists for each cube.
		\li Define the equivalent points and check points for each cube.
    \end{algleo}

    \linonumber \textsc{Step 2 Upward pass}
    \begin{algleo}
        \li \textbf{for} each leaf cube $C$ in \emph{postorder} 
			traversal of the tree \textbf{do}
        \begin{algleo}
			\li \textbf{if} $C$ is in the low frequency regime
			\begin{algleo}
				\li Compute the \emph{non-directional outgoing equivalent 
					densities} by Equation \eqref{eq_s2m} using the sources 
					inside $C$ (S2M).
			\end{algleo}
			\li \textbf{else} ($C$ is in the high frequency regime)
			\begin{algleo}
				\li Compute the \emph{directional outgoing equivalent 
					densities} by Equation \eqref{eq_s2m} for each directional 
					cone using the sources inside $C$ (S2M).
			\end{algleo}
			\li \textbf{end if}
        \end{algleo}
        \li \textbf{end for}
        \li \textbf{for} each non-leaf cube $C$ in \emph{postorder} 
			traversal of the tree \textbf{do}
        \begin{algleo}
			\li \textbf{if} $C$ is in the low frequency regime
			\begin{algleo}
				\li Compute the \emph{non-directional outgoing equivalent 
					densities} using the \emph{non-directional outgoing 
					equivalent densities} of its child cubes (M2M).
			\end{algleo}
			\li \textbf{else} ($C$ is in the high frequency regime)
			\begin{algleo}
				\li Compute the \emph{directional outgoing equivalent 
					densities} for each directional cone using 
					the equivalent densities of its child cubes (M2M).
			\end{algleo}
			\li \textbf{end if}
        \end{algleo}
        \li \textbf{end for}
    \end{algleo}

    \linonumber \textsc{Step 3 Downward pass}
    \begin{algleo}
        \li \textbf{for} each non-leaf cube $C$ in \emph{preorder} 
			traversal of the tree \textbf{do}
        \begin{algleo}
			\li \textbf{if} $C$ is in the low frequency regime
			\begin{algleo}
				\li Add to the \emph{downward check potentials} produced by
					the \emph{downward equivalent densities} in its 
					interaction list by Equation \eqref{eq_m2l} (M2L).
				\li Add to the \emph{downward check potentials} of its child 
					cubes (L2L).
			\end{algleo}
			\li \textbf{else} ($C$ is in the high frequency regime)
			\begin{algleo}
				\li Add to the \emph{directional incoming check potentials} 
					produced by the \emph{directional incoming equivalent 
					densities} in its interaction list by Equation 
					\eqref{eq_m2l} (M2L).
				\li Add to the \emph{directional incoming check potentials} 
					or the \emph{downward check potentials} of its child 
					cubes (L2L).
			\end{algleo}
			\li \textbf{end if}
        \end{algleo}
        \li \textbf{end for}
        \li \textbf{for} each leaf cube $C$ in \emph{preorder} traversal 
			of the tree \textbf{do}
        \begin{algleo}
            \li Evaluate the potentials on $C$ by Equation \eqref{eq_l2t} (L2T).
        \end{algleo}
        \li \textbf{end for}
    \end{algleo}

    \linonumber \textsc{Step 4 Near-field interaction}
    \begin{algleo}
        \li \textbf{for} each leaf cube $C$ in \emph{preorder} traversal 
			of the tree \textbf{do}
        \begin{algleo}
            \li Add to the potential the contribution of near field source 
			densities by Equation \eqref{eq_sum} (S2T).
        \end{algleo}
        \li \textbf{end for}
    \end{algleo}

\end{algleo}

\section{Further accelerating techniques}

The most time consuming step in the fast directional algorithm is the M2L
translation, since it has to be performed many times for each cube. 
Therefore, the accelerating technique for M2L can considerably improve 
the performance of the algorithm. Consider the M2L matrix in 
\eqref{eq_m2l}, although its numerically rank is $T(\varepsilon)$, the
number of \emph{directional outgoing equivalent points} and \emph{directional 
incoming check points} $O(T(\varepsilon))$ are often chosen to be much larger 
than $T(\varepsilon)$ in order to maintain the precision. It leads to that the 
dimension of the M2L matrices are fairly larger than their ranks, thus 
the M2L translations can be accelerated by low rank approximations.

In this section, a new accelerating technique similar with SArcmp in 
\cite{kifmbem, opm2l} is proposed, which accelerate the M2L translations
by first reducing the dimensions of all the matrices, then 
performing the low rank decomposition individually. 

\subsection{Matrix reduction for M2L}
\label{sbsec_matred}

Our new accelerating approach can be considered as an improved version for
SArcmp, since the main idea remains the same, and only the matrix reduction 
step is different. Therefore before presenting our new accelerating 
approach, the matrix reduction approach in the SArcmp 
for fast directional algorithm is introduced first.

\subsubsection{Matrix reduction in SArcmp}

In SArcmp\cite{kifmbem, opm2l}, to reduce the dimensions of M2L matrices 
$\bm{K}^{(1)}, \bm{K}^{(2)}, \cdots, \bm{K}^{(t)}$ in a directional cone, first 
they are collected in a row to form a ``fat'' matrix 
$$
\bm{K}_\text{fat} = \left[ \bm{K}^{(1)} \quad \bm{K}^{(2)} \quad 
\cdots \quad \bm{K}^{(t)} \right],
$$
and in a column to form a ``thin'' matrix 
$$
\bm{K}_\text{thin} = \left[ \bm{K}^{(1)}; \quad \bm{K}^{(2)}; \quad 
\cdots; \quad \bm{K}^{(t)} \right],
$$
where $t$ is the number of M2L matrices in a directional cone.
Then perform singular value decomposition (SVD) 
\begin{subequations}
\begin{equation}
	\bm{K}_\text{fat} = \bm{U} \bm{\Sigma} \left[ 
	{\bm{V}^{(1)}}^\text{H}, \quad {\bm{V}^{(2)}}^\text{H}, \quad 
	\cdots, \quad {\bm{V}^{(t)}}^\text{H} \right],
\end{equation}
\begin{equation}
	\bm{K}_\text{thin} = \left[ 
	{\bm{Q}^{(1)}}, \quad {\bm{Q}^{(2)}}, \quad 
	\cdots, \quad {\bm{Q}^{(t)}} \right] \bm{\Lambda} 
	\bm{R}^\text{H}.
\end{equation}
\end{subequations}

For each translating matrix $\bm{K}^{(i)}$, there is
\begin{equation}
	\bm{K}^{(i)} = \bm{U} \bm{U}^\text{H} \bm{K}^{(i)} \bm{R} \bm{R}^\text{H}.
\end{equation}
Truncate the columns in $\bm{U}$ and $\bm{R}$ corresponding to small 
singular values in $\bm{\Sigma}$ and $\bm{\Lambda}$, respectively, the 
equation becomes
\begin{equation} \label{eq_transm2l}
	\bm{K}^{(i)} \approx \tilde{\bm{U}} \tilde{\bm{U}}^\text{H} 
	\bm{K}^{(i)} \tilde{\bm{R}} \tilde{\bm{R}}^\text{H} 
	= \tilde{\bm{U}} \tilde{\bm{K}}^{(i)} \tilde{\bm{R}}^\text{H},
\end{equation}
where $\tilde{\bm{U}}$ and $\tilde{\bm{R}}$ are the compressing matrices, and
$\tilde{\bm{K}}^{(i)} =  \tilde{\bm{U}}^\text{H} \bm{K}^{(i)} \tilde{\bm{R}}$ 
is the reduced M2L matrix.

From the definition of the M2L matrices \eqref{eq_m2l} we know 
that $\bm{K}_\text{fat}$ can be viewed as the
evaluating matrix for potentials in $X$ produced by sources in $Y$, and 
$\bm{K}_\text{thin}$ can be viewed as the evaluating matrix for potentials in 
$Y$ produced by sources in $X$, where $X$ and $Y$ satisfies the directional 
parabolic separation condition, as shown in Figure \ref{fig_para_sep}. 
Therefore, $\bm{K}_\text{fat}$ and $\bm{K}_\text{thin}$ are also of rank 
$T(\varepsilon)$, and the dimension of the reduced M2L matrices 
$\tilde{\bm{K}}$ is $T(\varepsilon)$.

\begin{figure}[h]
	\centering
	\subfigure{
		\includegraphics[width=0.3\textwidth]{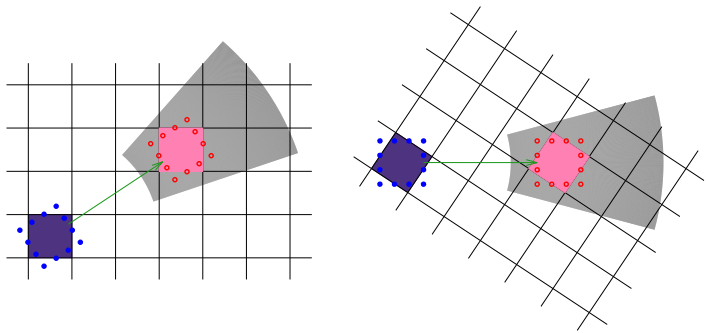}
	}
	\subfigure{
		\includegraphics[width=0.3\textwidth, origin=c, angle=-33.7]{figures/dirm2l}
	}
	\caption{M2L translation from $C$ (the purple cube) to $D$ (the pink cube) 
		in arbitrary directional cone.}
	\label{fig_dirm2l}
\end{figure}

To adapt SArcmp to fast directional algorithm, first let us consider 
two interacting cubes $C$ and $D$ in the high frequency regime. 
Assume $D$ is in an arbitrary $\bm{u}$-th directional cone of $C$, as 
illustrated in Figure \ref{fig_dirm2l}. The M2L matrix is 
also in the $\bm{u}$-th directional cone. Notice that 
the \emph{directional outgoing equivalent points} and the 
\emph{directional incoming check points} can be defined by 
rotating these points in the $(0, 0, 1)$-th directional cone.
Therefore the M2L matrix is the same with that in the $(0, 0, 1)$-th 
directional cone since the kernel is rotational invariant, as 
illustrated in Figure \ref{fig_dirm2l}. Therefore, all the 
M2L matrices in different directional cones are rotated to 
the $(0, 0, 1)$-th directional cone, and can be collected together
and compressed.

The key of the matrix reduction algorithm is finding the compressing
matrices $\tilde{\bm{U}}$ and $\tilde{\bm{R}}$. In SArcmp, all the M2L 
matrices has to be collected to generate $\tilde{\bm{U}}$ 
and $\tilde{\bm{R}}$. Assume the boundary is of size $D$, the size
of the cubes in the highest level is $w \sim \sqrt{D/k}$, thus
there are $O(D^2/w^2) = O(kD) = O(\sqrt{N})$ M2L matrices in the
highest level. Therefore, the computational complexity for collecting 
the M2L matrices and performing SVD is $O(N^{3/2})$, therefore this 
approach is not suitable for the fast directional algorithm.
A new compressing scheme is proposed in the following section, in 
which the compressing matrices can be generated in the upward and 
downward pass, thus the matrices collecting step can be omitted to 
avoid this drawback.

\subsubsection{Matrix reduction without matrices collecting}
\label{sbsbsec_matred}

In the multilevel fast algorithms, evaluating the potential $\bm{p}$ 
in a leaf cube induced by the source densities $\bm{q}$ in another leaf 
cube in its far field can be computed as
\begin{equation} \label{eq_trans}
	\bm{p = TLKMSq},
\end{equation}
where $\bm{S}, \bm{M}, \bm{K}, \bm{L}, \bm{T}$ are the S2M, M2M, M2L, 
L2L and L2T matrices, respectively. In the fast directional algorithm
base on equivalent densities, the M2M matrix $M$ and L2L matrix $L$ 
are also computed by two steps, which are similar with S2M and L2T, as 
shown in Figure \ref{fig_s2m} and \ref{fig_l2t},
\begin{subequations} \label{eq_matml}
\begin{equation}
	\bm{M} = \bm{G}_\text{up}^{+} \bm{E}_\text{up},
\end{equation}
\begin{equation}
	\bm{L} = \bm{E}_\text{dn} \bm{G}_\text{dn}^{+},
\end{equation}
\end{subequations}
where $\bm{G}_\text{up}$ is the matrix evaluating the \emph{directional outgoing 
check potentials} produced by the \emph{directional outgoing equivalent 
densities}, $\bm{G}_\text{dn}$ is the matrix evaluating the \emph{directional 
outgoing incoming potentials} produced by the \emph{directional incoming 
equivalent densities}, and the superscript ``+'' denotes the Moore-Penrose 
inverse. They can be computed via performing SVD for $\bm{G}_\text{up}$ and
$\bm{G}_\text{dn}$
\begin{subequations} \label{eq_svdg}
\begin{equation}
	\bm{G}_\text{up} = \bm{U}_\text{up} \bm{\Sigma}_\text{up}
	\bm{V}_\text{up}^\text{H}, 
\end{equation}
\begin{equation}
	\bm{G}_\text{dn} = \bm{U}_\text{dn} \bm{\Sigma}_\text{dn}
	\bm{V}_\text{dn}^\text{H}.
\end{equation}
\end{subequations}
Then the Moore-Penrose inverses $\bm{G}_\text{up}^{+}$ and 
$\bm{G}_\text{dn}^{+}$ can be approximated by inverting \eqref{eq_svdg} and
truncating the columns corresponding to tiny singular values 
\begin{equation} \label{eq_epsg}
	\sigma_i < \varepsilon \sigma_0, 
\end{equation}
where $\sigma_0$ is the largest singular value of $\bm{G}^{+}$. Thus 
\begin{subequations} \label{eq_invg}
\begin{equation}
	\bm{G}_\text{up}^{+} \approx \tilde{\bm{V}}_\text{up} 
	\tilde{\bm{\Sigma}}_\text{up}^{-1} \tilde{\bm{U}}_\text{up}^\text{H}, 
\end{equation}
\begin{equation}
	\bm{G}_\text{dn}^{+} \approx \tilde{\bm{V}}_\text{dn} 
	\tilde{\bm{\Sigma}}_\text{dn}^{-1} \tilde{\bm{U}}_\text{dn}^\text{H}.
\end{equation}
\end{subequations}
Substituting \eqref{eq_invg} to \eqref{eq_matml} and \eqref{eq_trans}, one gets
\begin{equation}
\begin{split}
	\bm{LKM} = & \bm{E}_\text{dn} \bm{G}_\text{dn}^{+} \bm{K} 
	\bm{G}_\text{up}^{+} \bm{E}_\text{up} \\
	\approx & \bm{E}_\text{dn} \tilde{\bm{V}}_\text{dn} 
	\tilde{\bm{\Sigma}}_\text{dn}^{-1} \tilde{\bm{U}}_\text{dn}^\text{H} \bm{K}
	\tilde{\bm{V}}_\text{up} \tilde{\bm{\Sigma}}_\text{up}^{-1} 
	\tilde{\bm{U}}_\text{up}^\text{H} \bm{E}_\text{up}.
\end{split}
\end{equation}

Since the directional equivalent points and the directional check points
also satisfies the directional parabolic separation condition, as illustrated in
Figure \ref{fig_slfda}, the number of columns in truncated matrices 
$\tilde{\bm{U}}_\text{dn}$ and $\tilde{\bm{V}}_\text{up}$
is also $T(\varepsilon)$. Therefore, take $\tilde{\bm{U}}_\text{dn}$ and 
$\tilde{\bm{V}}_\text{up}$ as the compressing matrices, the M2L matrices can 
also be compressed into more compact form 
\begin{equation} \label{eq_cmpm2l}
	\tilde{\bm{K}} = \tilde{\bm{U}}_\text{dn}^\text{H} \bm{K} 
	\tilde{\bm{V}}_\text{up}.
\end{equation}

\subsection{Low rank decomposition for reduced M2L matrices}

After the matrix reduction in Section \ref{sbsec_matred}, the resulting M2L
matrices are still of low rank\cite{opm2l, kifmbem}, this means the M2L translations 
can be further accelerated by performing low rank approximations for the M2L 
matrices individually. In this paper, this is also done in the same way as that in
\cite{kifmbem}, i.e., the low rank decomposition is achieved by SVD.

For each $k$ dimensional reduced M2L matrix $\tilde{\bm{K}}_{k \times k}$, 
perform SVD
\begin{equation}
	\tilde{\bm{K}} = \bm{U} \bm{S} \bm{Q}^\text{H}.
\end{equation}
Then truncate the columns corresponding to the tiny singular values
\begin{equation} \label{eq_epsk}
	\sigma_i < \varepsilon \sigma_0,
\end{equation}
where $\sigma_0$ is the largest singular value of $\bm{K}$. Thus
\begin{equation}
	\tilde{\bm{K}}_{k \times k} \approx \hat{\bm{U}}_{k \times r} 
	\hat{\bm{S}}_{r \times r} \hat{\bm{Q}}_{r \times k}^\text{H}
	= \hat{\bm{U}}_{k \times r} \hat{\bm{V}}_{r \times k},
\end{equation}
where $\hat{\bm{V}}_{r \times k} = \hat{\bm{S}}_{r \times r} 
\hat{\bm{Q}}_{r \times k}^\text{H}$. The translations can be more efficient
when $r < \frac{1}{2}k$.

\begin{figure}[h]
	\centering
	\includegraphics[width=0.4\textwidth]{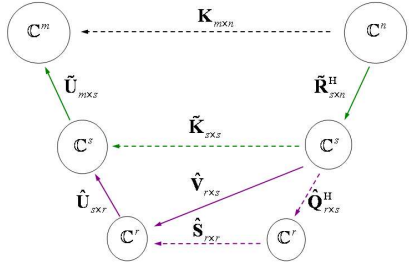}
	\caption{The accelerating approach for M2L.}
	\label{fig_sarcmp}
\end{figure}

The overall accelerating approach for M2L translations are illustrated in Figure
\ref{fig_sarcmp}. Where, $\bm{K}_{m \times n}$ is the original M2L 
matrix, $\tilde{\bm{U}}_{m \times k}$ and $\tilde{\bm{R}}_{k \times n}$
is the compressing matrices in Section \ref{sbsbsec_matred}, 
$\hat{\bm{U}}_{k \times r}$ and $\hat{\bm{V}}_{r \times k}$ are the
individually low rank decomposition matrices for M2L matrices.

\subsection{Accelerating technique for upward and downward passes}

It is proved in \cite{kifmbem} that, the compressing matrices for 
M2L reduction can also be used to compress the translation matrices in
upward and downward passes. In this paper, the compressing matrices 
are generated in computing the translation matrices in the upward and 
downward passes, and the compressing approach is simpler.

Assume $\bm{M}_l = \bm{G}_{\text{up}, l}^{+} \bm{E}_{\text{up}, l}$ is
the M2M matrix translating the \emph{directional outgoing equivalent 
densities} in the $(l+1)$-th level to the \emph{directional outgoing 
equivalent densities} in the $l$-th level, and $\bm{S}_{l+1} = 
\bm{G}_{\text{up}, l+1}^{+} \bm{E}_{\text{up}, l+1}$ is the S2M matrix
translating the sources inside a leaf cube in the $(l+1)$-th level to
the \emph{directional outgoing equivalent densities} in the $l$-th level.
Combining with Equations \eqref{eq_invg} and \eqref{eq_trans}, there is
\begin{equation}
\begin{split}
	\bm{M}_l \bm{S}_{l+1} = & \bm{G}_{\text{up}, l}^{+} \bm{E}_{\text{up}, l}
	\bm{G}_{\text{up}, l+1}^{+} \bm{E}_{\text{up}, l+1} \\
	\approx & \tilde{\bm{V}}_{\text{up}, l} \tilde{\bm{\Sigma}}_{\text{up}, l}^{-1}
	\tilde{\bm{U}}_{\text{up}, l}^\text{H} \bm{E}_{\text{up}, l}
	\tilde{\bm{V}}_{\text{up}, l+1} \tilde{\bm{\Sigma}}_{\text{up}, l+1}^{-1}
	\tilde{\bm{U}}_{\text{up}, l+1}^\text{H} \bm{E}_{\text{up}, l+1} \\
	= & \tilde{\bm{V}}_{\text{up}, l} (\tilde{\bm{\Sigma}}_{\text{up}, l}^{-1}
	\tilde{\bm{U}}_{\text{up}, l}^\text{H} \bm{E}_{\text{up}, l}
	\tilde{\bm{V}}_{\text{up}, l+1}) (\tilde{\bm{\Sigma}}_{\text{up}, l+1}^{-1}
	\tilde{\bm{U}}_{\text{up}, l+1}^\text{H} \bm{E}_{\text{up}, l+1}) \\
	= & \tilde{\bm{V}}_{\text{up}, l} \tilde{\bm{M}}_{\text{up}, l}
	\tilde{\bm{S}}_{\text{up}, l+1}.
\end{split}
\end{equation}
Where, $\tilde{\bm{V}}_{\text{up}, l}$ would be used to compress the M2L matrices
in $l$-th level, as shown in \eqref{eq_cmpm2l},
\begin{equation}
	\tilde{\bm{M}}_{\text{up}, l} = \tilde{\bm{\Sigma}}_{\text{up}, l}^{-1}
	\tilde{\bm{U}}_{\text{up}, l}^\text{H} \bm{E}_{\text{up}, l}
	\tilde{\bm{V}}_{\text{up}, l+1}
\end{equation}
is the compressed M2M matrix in $l$-th level, and
\begin{equation}
	\tilde{\bm{S}}_{\text{up}, l+1} = \tilde{\bm{\Sigma}}_{\text{up}, l+1}^{-1}
	\tilde{\bm{U}}_{\text{up}, l+1}^\text{H} \bm{E}_{\text{up}, l+1}
\end{equation}
is the compressed S2M matrix in $(l+1)$-th level. The translating
matrices in the downward pass compressed in the same manner, resulting
\begin{equation}
	\tilde{\bm{L}}_{\text{dn}, l} = 
	\tilde{\bm{U}}_{\text{dn}, l+1}^\text{H} \bm{E}_{\text{dn}, l}
	\tilde{\bm{V}}_{\text{dn}, l} \tilde{\bm{\Sigma}}_{\text{dn}, l}^{-1}
\end{equation}
be the compressed L2L matrix translating the \emph{directional incoming
check potentials} in the $l$-th level to the \emph{directional incoming
check potentials} in the $(l+1)$-th level, and 
\begin{equation}
	\tilde{\bm{T}}_{\text{dn}, l+1} = \bm{E}_{\text{dn}, l+1}
	\tilde{\bm{V}}_{\text{dn}, l+1} \tilde{\bm{\Sigma}}_{\text{dn}, l+1}^{-1}
\end{equation}
be the compressed L2T matrix translating the \emph{directional incoming
check potentials} of a leaf cube in the $(l+1)$-th level to the target potentials.

\section{Numerical studies}

The performance of our fast directional algorithm for Burton-Miller 
formulation is demonstrated by several numerical examples. 
The codes are implemented serially in C++.
The Burton-Miller formulation \eqref{eq_bm}
is discretized by collocation method and piecewise constant elements, i.e., %
assume $\bm{x}_i$ is the centroid of the $i$-th triangular element 
$\triangle_i$, the weight functions $w_i(\bm{x})$ and the basis 
functions $\chi_i(\bm{x})$ are chosen to be 
\begin{equation}
\begin{split}
	w_i(\bm{x}) =& \delta(\bm{x}_i), \\
	\chi_i(\bm{x}) =& \delta_{ij}, \quad \bm{x} \in \triangle_j.
\end{split}
\end{equation}
The resulting linear systems are solved by GMRES solver, and its
converging tolerance is set to be equal to the singular value 
truncating threshold $\varepsilon$ in \eqref{eq_epsg} and 
\eqref{eq_epsk}. All the computational results are computed on 
a computer with a Xeon 5450 (2.66 GHz) CPU and 32 GB RAM.

\subsection{Performance of our algorithm}

First let us study the performance of our algorithm, in which the 
equivalent points and the check points are distributed 
straightforwardly, and the translation operators are compressed. 
It is studied by a unit sphere pulsating problem. That is, the 
surface of the unit sphere pulsates with uniform radial velocity 
$v_a = \frac{\partial u}{\partial \bm{n}} = 1$. The numerical error 
of the velocity potential on the surface can be computed via the 
analytical solution
\begin{equation}
	u = \frac{1}{1 + ik}.
\end{equation}

The unit sphere is first discretized into $N=512$ triangular elements and it is
used to compute the pulsating problem with $k=\pi$. That is, the 
diameter equals 1 wavelength. Then the mesh is refined and the wavenumber
is doubled 6 times. The finest mesh has $N=2097152$ elements and the 
diameter of the sphere is 64 in terms of wavelength.

First let us study the influence of the number of equivalent points along 
each direction $p$ on the accuracy and efficiency of the algorithm.
Let $p = \log(1/\varepsilon) + p_0$. The performance of the 
algorithm for different $p_0$
is studied by choosing $p_0=0, 1, 2, 3$ for the unit sphere pulsating 
problem with the finest mesh and $\varepsilon = \text{1e-3}$. The 
results are listed in Table \ref{tab_varp0}, where $T_\text{t}$
is the total time cost, $T_\text{it}$ is the running time 
in each iteration, $N_\text{it}$ is the number of iterations, and 
$M$ is the memory consume. It is shown 
that the resulting error maintains almost the
same when $p_0 \ge 1$, while the total time cost is considerably
increased. Therefore, we choose $p_0=1$ for the following examples.

\begin{table} [ht]
	\centering
	\caption{Results of the unit sphere pulsating problem with different
	   number of equivalent points $N=2097152, \varepsilon = \text{1e-3}$.}
	\vspace{1em}
	\begin{tabular}{rrrrrr}
		\hline
		$p_0$     & $T_\text{t}(s)$	& $T_\text{it}$(s)	& $N_\text{it}$	& $M$(MB)	& $L_2$-error \\
		\hline
		0	& 12067.8	& 274.75	& 8	& 24451.2	& 1.90e-2\\
		1	& 14250.9	& 260.45	& 4	& 27822.7	& 3.90e-3\\
		2	& 18268.9	& 260.20	& 4	& 27846.2	& 3.37e-3\\
		3	& 23543.5	& 325.01	& 4	& 27526.2	& 3.24e-3\\
		\hline
	\end{tabular}
	\label{tab_varp0}
\end{table}

The results for the unit sphere pulsating problem are listed in Table 
\ref{tab_usphpls}. The case with $N=2097152, \varepsilon=\text{1e-4}$
is not computed because of memory constraint. It is shown that 
the complexity of our algorithm grows 
almost linearly with respect to the number of degrees. The resulting
error concludes the discretization error of BEM and the approximating 
error $\varepsilon$ of the fast algorithm. The resulting error for 
meshes with less than $32768$ elements preserves almost the same with
$\varepsilon = \text{1e-3}$ and 1e-4, which shows that the resulting error
is bounded by the precision of the boundary integral discretization. 
The errors for finer meshes maintain almost the same with 
$\varepsilon=\text{1e-3}$, this is because they are bounded by the 
accuracy of the fast accelerating scheme. The error continue 
decreasing at the same rate with $\varepsilon=\text{1e-4}$, which indicates
that the error of the fast accelerating scheme is smaller and more precise
results can be obtained by decreasing $\varepsilon$. That is, the error 
can be reduced linearly to $O(\varepsilon)$, which indicates that
our fast algorithm is quite stable.

\begin{table} [ht]
	\centering
	\caption{Results of the unit sphere pulsating problem}
	\vspace{1em}
	\begin{tabular}{rrrrrrrr}
		\hline
		$N$     & $k$	& $T_\text{t}$(s)	& $T_\text{it}$(s)	& $N_\text{it}$	& $M$(MB)	& $L_2$-error \\
		\hline
		\multicolumn{7}{c}{$\varepsilon=$1e-3}\\
		512		& $\pi$		&     1.76	&   0.02	& 3	&     4.59	& 3.65e-2\\
		2048	& $2\pi$	&     4.26	&   0.07	& 3	&    19.99	& 1.86e-2\\
		8192	& $4\pi$	&    50.95	&   0.46	& 3	&   119.76	& 9.26e-3\\
		32768	& $8\pi$	&   166.85	&   2.63	& 3	&   346.67	& 4.84e-3\\
		131072	& $16\pi$	&   962.98	&  13.79	& 3	&  2088.75	& 3.11e-3\\
		524288	& $32\pi$	&  3427.14	&  69.15	& 3	&  6761.24	& 2.92e-3\\
		2097152	& $64\pi$	& 14250.90	& 260.45	& 4	& 27822.70	& 3.90e-3\\
		\multicolumn{7}{c}{$\varepsilon=$1e-4}\\                  
		512		& $\pi$		&     2.63	&   0.02	& 6	&     8.98	& 3.65e-2\\
		2048	& $2\pi$	&     8.53	&   0.16	& 5	&    39.66	& 1.85e-2\\
		8192	& $4\pi$	&    54.89	&   0.74	& 4	&   192.34	& 9.27e-3\\
		32768	& $8\pi$	&   254.08	&   5.66	& 4	&   499.10	& 4.63e-3\\
		131072	& $16\pi$	&  1490.30	&  29.58	& 5	&  2920.10	& 2.35e-3\\
		524288	& $32\pi$	&  5457.84	& 143.36	& 6	&  9223.96	& 1.32e-3\\
		\hline
	\end{tabular}
	\label{tab_usphpls}
\end{table}

\subsection{Comparison to the wideband FMM}

To compare the performance of the current fast directional 
algorithm with the wideband fast multipole method, the unit 
sphere scattering problem in Section 4.1 in \cite{wbfmbem} 
is computed. The point source is at (-2, 0, 0). Two cases with
$k = 5.0$ and $k = 50$ are computed. The sphere surface is 
discretized into approximately the same number of elements with 
that in \cite{wbfmbem}, and we chose $\varepsilon = 1e-3$. The 
results are listed in Table \ref{tab_cmpwbfmbem}.
Note that the wideband FMM is parallelized and performed by a four core 
computer with a 64-bit Intel Core$^\text{TM}$ 2 Duo CPU, thus it should 
be about 4 times faster. However, it is shown that, our algorithm consumes
almost the same time when $ka=50.0$ and is much faster when $ka=5.0$. 
Therefore, our algorithm is much more efficient than the wideband 
fast multipole algorithm.

\begin{table} [ht]
	\centering
	\caption{Comparison of results for sphere scattering problems reported
	in Ref. \cite{wbfmbem} and the present study.}
	\vspace{1em}
	\begin{tabular}{crrcr}
		\hline
		~	& $k$		& $N$		& cores	& $T_\text{it}$(s)\\
		\hline
		\cite{wbfmbem}	& 5.0	& 101270	& 4	&  7.08	\\
		Present			& 5.0	& 106032	& 1	&  2.15	\\
		\cite{wbfmbem}	& 50.0	& 101270	& 4	& 10.84	\\
		Present			& 50.0	& 106032	& 1	& 12.38	\\
		\hline
	\end{tabular}
	\label{tab_cmpwbfmbem}
\end{table}

\subsection{Plane wave scattering problems}

For plane wave scattering problems, two sound-hard obstacles are 
considered, and the incident plane wave is assumed to be propagating 
in (1, 0, 0) direction. The singular truncating threshold and the 
GMRES converging tolerance is set to be $\varepsilon = \text{1e-3}$.
The triangular meshes used in the examples have a refinement of approximately
16 elements per wavelength. 

\begin{figure}[h]
	\centering
	\subfigure[Acoustic field.]{
		\includegraphics[height=0.4\textwidth]{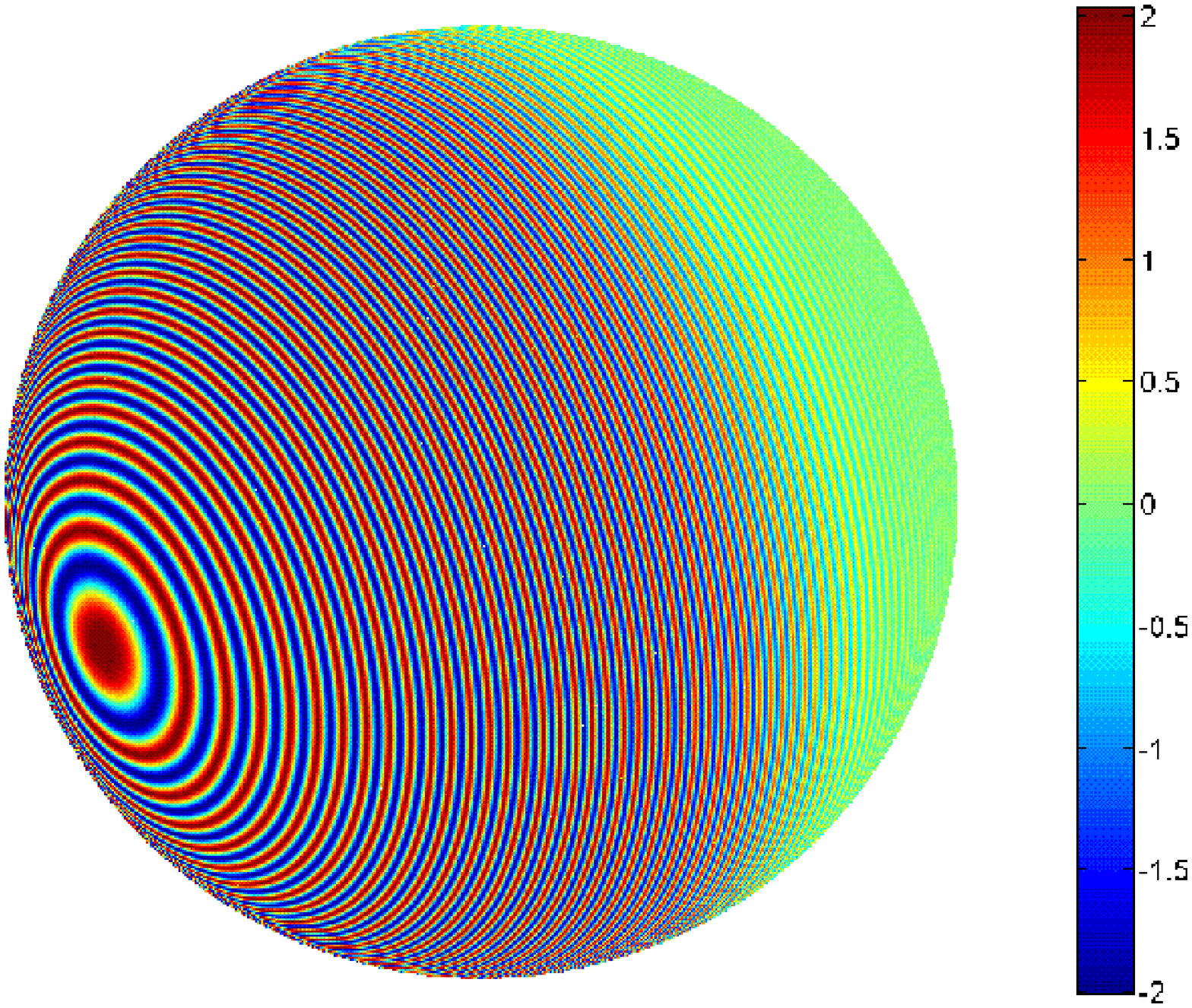}
		\label{fig_usphfield}
	}
	\subfigure[Scattering velocity potential.]{
		\includegraphics[height=0.4\textwidth]{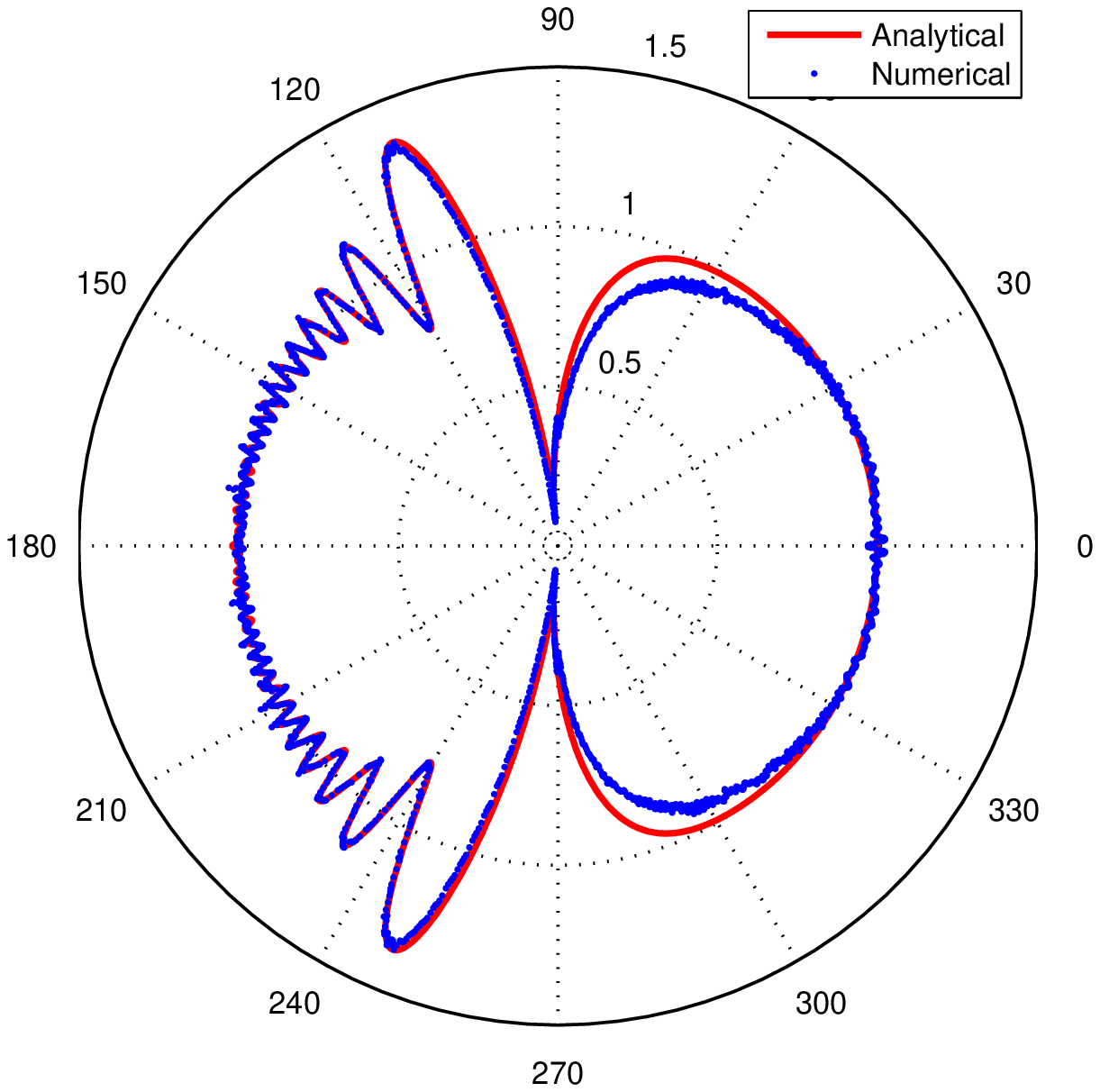}
		\label{fig_usphsca}
	}
	\caption{Resulting of the unit sphere scattering problem.}
	\label{fig_usphscares}
\end{figure}

First the sound-hard unit sphere scattering problem is calculated.
The sphere diameter is 64 wavelengths, i.e., $k = 201$. The sphere
surface is discretized into 2097152 triangular elements. The total time
cost for solving the problem $T_\text{t} =$ 8690s, the time cost 
in each iteration $T_\text{it} =$ 260.3s, the number of iterations 
$N_\text{it} = 14$, the memory consumption 
$M = $ 26.2GB. The resulting acoustic velocity potential on the 
surface is illustrated in Fig. \ref{fig_usphfield}, and the scattering 
field on the surface is illustrated in Fig. \ref{fig_usphsca}. 
It is shown that our numerical results agrees quite well with 
the analytical solution.

\begin{figure}[h]
	\centering
	\includegraphics[width=0.7\textwidth]{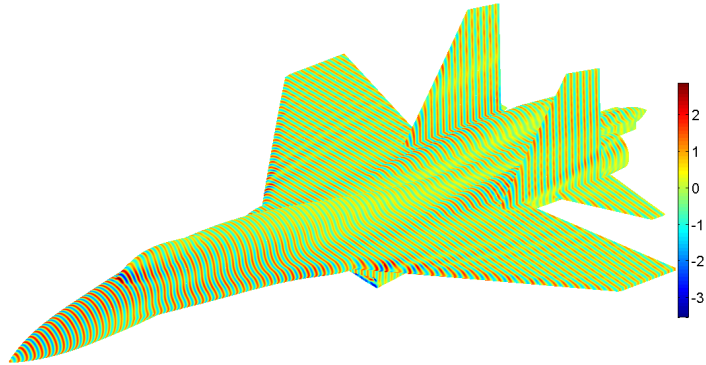}
	\caption{Resulting of the Su-27 scattering problem.}
	\label{fig_su27}
\end{figure}

The second model is a Su-27 fighter which is 21.49 meters long. 
We set $k = 30$, thus its size is 102.7 wavelengths. The surfaces 
are discretized into 4143908 triangular elements. The total time 
cost for solving the problem $T_\text{t} =$ 38799.3s, the time 
cost in each iteration $T_\text{it} =$ 133.57s, the number of iterations
$N_\text{it} = 256$ without any preconditioner.
The memory consumption $M = $ 25.08GB. The resulting acoustic velocity 
potential on the surface is illustrated in Fig. \ref{fig_su27}.

\begin{figure}[h]
	\centering
	\includegraphics[width=0.7\textwidth]{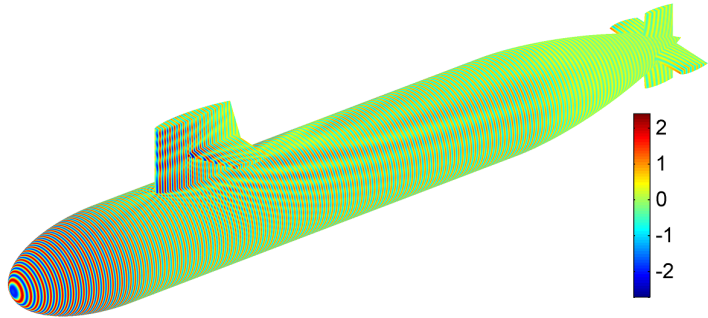}
	\caption{Resulting of a submarine scattering problem.}
	\label{fig_submarine}
\end{figure}

The second model is a submarine which is 80 meters long. 
We set $k = 12.5$, thus $kD = 1000$ and  its size is 159 wavelengths. 
The surfaces are discretized into 4041088 triangular elements. 
The total time cost for solving the problem $T_\text{t} =$ 13648s $\approx$ 
3.8h, the time 
cost in each iteration $T_\text{it} =$ 204.14s, the number of iterations
$N_\text{it} = 35$ without any preconditioner.
The memory consumption $M = $ 24.73GB. The resulting acoustic velocity 
potential on the surface is illustrated in Fig. \ref{fig_submarine}.

\section{Conclusion}

In this paper, the fast directional algorithm is adapted to accelerate
the acoustic problem computations with Burton-Miller formulation. 
Although there are four integrals in the Burton-Miller formulation, they 
can be evaluated efficiently by two fast summing approach. The outgoing
equivalent points and the outgoing check points are sampled directly instead 
of by pseudo skeleton approach, resulting in a simpler fast directional algorithm.
Then all the translations are accelerated by matrix reduction and low-rank 
approximation, which is similar with the SArcmp approach, while the 
compressing matrices are generated in advance, and no matrix collection
is required. The accuracy and efficiency of 
the algorithm are examined by numerical results. It is shown that the acoustic
scattering problems with over 4 million DOFs and $kD=1000$ can be computed in 
less than 4 hours.

\section*{Acknowledgements}

This work was supported by the Doctorate Foundation of Northwestern
Polytechnical University under Grant No. CX201220, National Science
Foundations of China under Grants 11074201 and 11102154, and Funds 
for Doctor Station from the Chinese Ministry of
Education under Grants 20106102120009 and 20116102110006.

\bibliographystyle{unsrt}
\bibliography{fda}

\begin{thebibliography}{10}

\bibitem{fmm}
L.~Greengard and V.~Rokhlin.
\newblock A fast algorithm for particle simulation.
\newblock {\em Journal of Computational Physics}, 73:325--348, 1987.

\bibitem{pwfmm}
L.~Greengard and V.~Rokhlin.
\newblock {A new version of the fast multipole method for the Laplace equation
  in three dimensions}.
\newblock {\em Acta Numerica}, 6(1):229--269, 1997.

\bibitem{fmm_Darve2000}
Eric Darve.
\newblock The fast multipole method: Numerical implementation.
\newblock {\em Journal of Computational Physics}, 160:195--240, 2000.

\bibitem{Hmatrix}
W.~Hackbusch.
\newblock {A sparse matrix arithmetic based on H-matrices. Part I: Introduction
  to H-matrices}.
\newblock {\em Computing}, 62:89--108, 1999.

\bibitem{pnlclst}
W.~Hackbusch and Z.~P. Nowak.
\newblock {On the fast matrix multiplication in the boundary element method by
  panel clustering}.
\newblock {\em Numerische Mathematik}, 54:463--491, 1989.

\bibitem{aca}
Mario Bebendorf.
\newblock Approximation of boundary element matrices.
\newblock {\em Numerische Mathematik}, 86:565--589, 2000.

\bibitem{aca_BR2003}
Mario Bebendorf and Sergej Rjasanow.
\newblock Adaptive low-rank approximation of collocation matrices.
\newblock {\em Computing}, 70:1--24, 2003.

\bibitem{wbem}
G.~Beylkin, R.~Coifman, and V.~Rokhlin.
\newblock Fast wavelet transforms and numerical algorithms.
\newblock {\em Pure Appl. Math.}, 37:141--183, 1991.

\bibitem{twbem}
Johannes Tausch.
\newblock A variable order wavelet method for the sparse representation of
  layer potentials in the non-standard form.
\newblock {\em Journal of Numerical Mathematics}, 12(3):233--254, 2004.

\bibitem{pfft}
Joel~R. Phillips and Jocob~K. White.
\newblock {A precorrected-FFT method for electrostatic analysis of complicated
  3-D structures}.
\newblock {\em IEEE Transactions on Computer-Aided Design of Integrated
  Circuits and Systems}, 16(10):1059--1072, 1997.

\bibitem{pfft_YanzyGaoxw}
Zai~You Yan and Xiao~Wei Gao.
\newblock {The development of the pFFT accelerated BEM for 3-D acoustic
  scattering problems based on the Burton and Miller's integral formulation}.
\newblock {\em Engineering Analysis with Boundary Elements}, 37:409--418, 2013.

\bibitem{pfft_Bruno2001}
Oscar~P. Bruno and Leonid~A. Kunyansky.
\newblock A fast, high-order algorithm for the solution of surface scattering
  problems: Basic implementation, tests, and applications.
\newblock {\em Journal of Computational Physics}, 110:80--110, 2001.

\bibitem{pfft_Bruno2012}
Oscar~P. Bruno, Tim Elling, and Catalin Turc.
\newblock Regularized integral equations and fast high-order solvers for
  sound-hard acoustic scattering problems.
\newblock {\em International Journal of Numerical Methods in Engineering},
  91:1045--1072, 2012.

\bibitem{diagonalform}
V.~Rokhlin.
\newblock {Diagonal forms of translation operators for the Helmholtz equation
  in three dimensions}.
\newblock {\em Applied and Computational Harmonic Analysis}, 1:82--93, 1993.

\bibitem{wbfmm}
Hongwei Cheng, William~Y. Crutchfield, Zydrunas Gimbutas, Leslie~F. Greengard,
  J.~Frank Ethridge, Jingfang Huang, Vladimir Rokhlin, Norman Yarvin, and
  Junsheng Zhao.
\newblock {A wideband fast multipole method for the Helmholtz equation in three
  dimensions}.
\newblock {\em Journal of Computational Physics}, 216:300--325, 2006.

\bibitem{wbfmm_Gumerov2009}
Nail~A. Gumerov and Ramani Duraiswami.
\newblock {A broadband fast multipole accelerated boundary element method for
  the three dimensional Helmholtz equation}.
\newblock {\em Journal of the Acoustical Society of America}, 125(1):191--205,
  2009.

\bibitem{fda}
Bj\"{o}rn Engquist and Lexing Ying.
\newblock Fast directional multilevel algorithms for oscillatory kernels.
\newblock {\em SIAM J. Sci. Comput.}, 29(4):1710--1737, 2007.

\bibitem{fda2010}
Bj\"{o}rn Engquist and Lexing Ying.
\newblock {Fast directional algorithms for the Helmholtz kernel}.
\newblock {\em Journal of Computational and Applied Mathematics},
  234:1851--1859, 2010.

\bibitem{fda_Chebyshev}
Matthias Messner, Martin Schanz, and Eric Darve.
\newblock {Fast directional multilevel summation for oscillatory kernels based
  on Chebyshev interpolation}.
\newblock {\em Journal of Computational Physics}, 231(4):1175--1196, 2012.

\bibitem{dirH2}
Mario Bebendorf, Christian Kuske, and Raoul Venn.
\newblock {Wideband nested cross approximation for Helmholtz problems}.
\newblock Number 536, Bonn, December 2012.

\bibitem{condBM}
R.~Kress.
\newblock Minimizing the condition number of boundary integral operators in
  acoustic and electromagnetic scattering.
\newblock {\em Quarterly Journal of Mechanics \& Applied Mathematics},
  38(2):323--341, 1985.

\bibitem{kifmm}
Lexing Ying, George Biros, and Denis Zorin.
\newblock A kernel-independent adaptive fast multipole algorithm in two and
  three dimensions.
\newblock {\em Journal of Computational Physics}, 196:591--626, 2004.

\bibitem{kifmbem}
Yanchuang Cao, Lihua Wen, and Junjie Rong.
\newblock {A SVD accelerated kernel-independent fast multipole method and its
  application to BEM}.
\newblock {\em arXiv: 1211.2517v2}, pages 1--19, 2012.

\bibitem{opm2l}
Matthias Messner, Berenger Bramas, Olivier Coulaud, and Eric Darve.
\newblock {Optimized M2L kernels for the Chebyshev interpolation based fast
  multipole method}.
\newblock {\em arXiv preprint, arXiv:1210.7292}, pages 1--23, 2012.

\bibitem{wbfmbem}
W.~R. Wolf and S.~K. Lele.
\newblock Wideband fast multipole boundary element method: Application to
  acoustic scattering from aerodynamic bodies.
\newblock {\em International journal for numerical methods in fluids},
  67:2108--2129, 2011.

\end{thebibliography}

\end{document}